\documentclass[12pt]{article}
\usepackage{amscd,amsfonts,amssymb,amsmath,latexsym,array,hhline}
\usepackage[dvips]{graphics}
\makeatletter
\@addtoreset{equation}{section}
\renewcommand{\@begintheorem}[2]{\begin{trivlist}\it
\item[\hspace{\labelsep}{\bf #1\ #2.}]}
\renewcommand{\@opargbegintheorem}[3]{\begin{trivlist}\it
\item[\hspace{\labelsep}{\bf #1\ #2\ (#3).}]}
\renewcommand{\@endtheorem}{\end{trivlist}}
\renewcommand{\@cite}[2]{[{#1\if@tempswa ; #2\fi}]}
\makeatother

\newcommand{\paragr}{\hspace{6mm}}

\newcommand{\skm}{\bigskip}
\newcommand{\skb}{\bigskip}

\newcommand{\al}{\alpha}
\newcommand{\be}{\beta}
\newcommand{\ga}{\gamma}

\newcommand{\De}{\Delta}

\newcommand{\la}{\lambda}
\newcommand{\si}{\sigma}

\newcommand{\eps}{\varepsilon}
\renewcommand{\phi}{\varphi}
\renewcommand{\kappa}{\varkappa}
\newcommand{\N}{\mathbb{N}}

\newcommand{\R}{\mathbb{R}}

\newcommand{\F}{\mathcal{F}}

\newcommand{\CCC}{\mathcal{C}}
\newcommand{\DDD}{\mathcal{D}}

\newcommand{\RRR}{\mathcal{R}}

\newcommand{\EX}{\mathcal{X}}

\newcommand{\EE}{\mathsf{E}}
\newcommand{\PP}{\mathsf{P}}
\newcommand{\QQ}{\mathsf{Q}}
\newcommand{\RR}{\mathsf{R}}

\newcommand{\emp}{\emptyset}
\newcommand{\lb}{\langle}
\newcommand{\rb}{\rangle}

\newcommand{\wl}{\overline}
\newcommand{\xra}{\xrightarrow}
\newcommand{\da}{\downarrow}

\newcommand{\ds}{\displaystyle}
\newcommand{\cond}{\hspace{0.3mm}|\hspace{0.3mm}}

\newcommand{\argmin}{\mathop{\rm argmin}}

\newcommand{\cl}{\mathop{\rm cl}}

\renewcommand{\inf}{\mathop{\rm inf\rule[-0.8mm]{0mm}{1mm}}}

\newenvironment{mitemize}%
{\begin{list}{$\bullet$}{
\leftmargin=32pt
\rightmargin=0pt
\labelsep=5pt
\labelwidth=20pt
\itemindent=0pt
\topsep=5pt plus 2pt minus 4pt
\partopsep=2pt plus 1pt minus 1pt
\parsep=0pt
\itemsep=0pt}}%
{\end{list}}

\begin{document}
\vspace*{10mm}
\begin{center}\bf
CAPM, REWARDS, AND EMPIRICAL ASSET PRICING

\vspace{2mm}
WITH COHERENT RISK
\end{center}

\begin{center}\itshape\bfseries
Alexander S.~Cherny$^*$,\quad Dilip B.~Madan$^{**}$
\end{center}

\begin{center}
\textit{$^*$Moscow State University}\\
\textit{Faculty of Mechanics and Mathematics}\\
\textit{Department of Probability Theory}\\
\textit{119992 Moscow, Russia}\\
\texttt{E-mail: cherny@mech.math.msu.su}\\
\texttt{Webpage: http://mech.math.msu.su/\~{}cherny}
\end{center}

\begin{center}
\textit{$^{**}$Robert H.~Smith School of Business}\\
\textit{Van Munching Hall}\\
\textit{University of Maryland}\\
\textit{College Park, MD 20742}\\
\texttt{E-mail: dmadan@rhsmith.umd.edu}\\
\texttt{Webpage: http://www.rhsmith.umd.edu/faculty/dmadan}
\end{center}

\begin{abstract}
\textbf{Abstract.}
The paper has 2 main goals:
\begin{mitemize}
\item[\bf 1.] We propose a variant of the CAPM based on
coherent risk.
\item[\bf 2.] In addition to the real-world measure
and the risk-neutral measure, we propose the third
one: \textit{the extreme measure}.
The introduction of this measure provides a powerful
tool for investigating the relation between the first
two measures. In particular, this gives us
\begin{mitemize}
\item a new way of measuring reward;
\item a new approach to the empirical asset pricing.
\end{mitemize}
\end{mitemize}

\end{abstract}

\section*{Contents}

\noindent\hbox to \textwidth{\ref{I}\ \ Introduction
\dotfill \pageref{I}}

\noindent\hbox to \textwidth{\ref{CRM}\ \ Coherent
Risk Measures \dotfill \pageref{CRM}}

\noindent\hbox to \textwidth{\ref{CAPM}\ \ CAPM
\dotfill \pageref{CAPM}}

\noindent\hbox to \textwidth{\ref{REAP}\ \ Rewards and
Empirical Asset Pricing \dotfill \pageref{REAP}}

\noindent\hbox to \textwidth{\ref{SC}\ \ Summary and
Conclusion \dotfill \pageref{SC}}

\noindent\hbox to \textwidth{References
\dotfill \pageref{R}}

\section{Introduction}
\label{I}

\textbf{1. CAPM.}
The Capital Asset Pricing Model is
based on variance as the measure of risk.
However, it has been clear from the outset that this
way of measuring risk has a serious drawback: it penalizes
high profits in exactly the same way as high losses.
Already in his 1959 book~\cite{M59} Markowitz suggested
semivariance as a substitute for variance.
But although semivariance is wiser than variance as a
measure of risk, it is less convenient analytically
and for this reason did not find its way to finance.

In 1997, a fundamentally new way of measuring risk
was proposed by Artzner, Delbaen, Eber, and Heath~\cite{ADEH97},
\cite{ADEH99}. They introduced the concept of a
\textit{coherent risk measure}.
In our opinion, these measures of risk are definitely
wiser than the standard ones, i.e. variance and V@R.
These new risk measures have already occupied a considerable
part of the modern financial mathematics (see the
literature review in~\cite{CM061}), and their theory is
progressing at an impressive speed.
We believe that within a few decades they will occupy a
firm position in practice also.
The theory of coherent risk measures is already termed
in some sources the ``third revolution in finance''
(see~\cite{Sz04}).

A coherent risk measure is a function on the random
variables of the form
\begin{equation}
\label{i1}
\rho(X)=-\inf_{\QQ\in\DDD}\EE_\QQ X,
\end{equation}
where $\DDD$ is a set of probability measures termed
\textit{probabilistic scenarios}.
From the financial point of view, $X$ is the
discounted P\&L produced by some portfolio over the unit time period.
Thus, $\rho(X)$ is the minimal capital needed
for the expected terminal wealth to be positive under
each scenario.

The above definition is very general.
For practical purposes one needs to select a convenient
subclass.
One of the best subclasses of coherent risk measures
known so far is \textit{Tail V@R} defined as follows:
$$
\rho_\la(X)=-\EE(X\cond X\le q_\la),
$$
where $\la\in(0,1]$ is a fixed number and $q_\la$ is
the $\la$-quantile of~$X$ (one can check that this is
indeed a coherent risk measure).
However, its empirical estimation might be problematic
due to the scarcity of tail events.
In~\cite{CM061}, we proposed a coherent risk measure,
which has a very clear meaning and admits a very simple
estimation procedure. We called it \textit{Alpha V@R}.
It is a risk measure of the form
$$
\rho_\al(X)=-\EE\min_{i=1,\dots,\al}X_i,
$$
where $\al\in\N$ is a fixed number and $X_1,\dots,X_\al$
are independent copies of~$X$.
We believe that the family of Alpha V@Rs (indexed by~$\al$)
is the best one-parameter family of coherent risk measures.
We also introduced in~\cite{CM061} the class
\textit{Beta V@R}, which, in our opinion,
is the best two-parameter family of coherent risks.

The first basic task of this paper is:
\textit{To build the CAPM based on coherent risk measures}.
The three main results of the CAPM are:
\begin{mitemize}
\item Establishing the SML relation.
\item Finding the form of an agent's optimal portfolio.
\item Finding the market risk premium.
\end{mitemize}
We provide coherent-based counterparts of these results.

\skb
\textbf{2. Rewards and empirical asset pricing.}
A very hot topic in the modern academic finance literature
is the relationship between the real-world and the
risk-neutral measures.
One trend consists in estimating the real-world measure
from the risk-neutral one (see~\cite{BKM03}, \cite{BP04},
\cite{LSTX05}).
Another trend is comparing the real-world and the risk-neutral
measures to derive the empirical pricing kernel
(see~\cite{AL00}, \cite{J00}, \cite{RE02}, \cite{RJ04}).

Our results on the CAPM show that, in addition to the
real-world measure and the risk-neutral measure, the
third measure plays a fundamental role:
$$
\text{\itshape\bfseries the extreme measure.}
$$

The notion of an extreme measure was introduced in~\cite{C061}
in order to define the coherent risk contribution.
Consider a firm measuring its risk by a
coherent risk measure~$\rho$ and producing a discounted
P\&L~$W$ over the unit time period.
The corresponding extreme measure is defined as
$$
\QQ(W)=\argmin_{\QQ\in\DDD}\EE_\QQ W,
$$
where $\DDD$ is the set standing in~\eqref{i1}.
Let $X$ be the discounted P\&L produced by some trade
over the same period. As shown in~\cite[Subsect.~2.5]{C061},
if $X$ is small as compared to~$W$, then
$$
\rho(W+X)-\rho(W)\approx-\EE_{\QQ(W)}X.
$$
It is seen from this relation that the notion of an
extreme measure is very useful for the risk measurement
purposes.
As found in~\cite{CM061}, this notion is also very
useful for the risk management purposes:
we prove that if the desks of
a firm are measuring coherent risk contributions rather
than outstanding risks and the desks are allowed to trade
risk limits between them, then the desks find
themselves the globally optimal portfolio.
Furthermore, the fair price intervals corresponding to
various pricing techniques (optimality pricing,
equilibrium pricing) considered in~\cite{C062} are expressed
through extreme measures;
the solution of the capital allocation problem considered
in~\cite{C061} is expressed through extreme measures;
see also our paper~\cite{CM062}, where the extreme measure
serves as one of the most natural examples of valuation measures.

In this paper, we obtain the following relationship:
\begin{equation}
\label{i2}
\RR=\frac{1}{1+R_*}\,\PP+\frac{R_*}{1+R_*}\,\QQ.
\end{equation}
Here $\PP$ is the real-world measure,
$\QQ$ is the market extreme measure
(i.e. this is the extreme measure corresponding to the
change of the market index like S\&P~500),
$\RR$ is the risk-neutral measure,
and $R_*$ is the risk premium for the market portfolio.
To be more precise, $\RR$ is a particular representative
of the class of risk-neutral measures, which we term
the \textit{contact measure}. It is closely connected
with the \textit{No Better Choice pricing} introduced
in~\cite[Subsect.~3.1]{C062}.
The relationship between $\PP$, $\QQ$, and $\RR$ is
illustrated by the following diagram:

\begin{picture}(150,12)(-14,-4)
\put(0,0){\circle*{1}}
\put(5,0){\circle*{1}}
\put(120,0){\circle*{1}}
\put(0,0){\line(1,0){120}}
\put(-1.5,2){$\PP$}
\put(3.5,2){$\RR$}
\put(118,2){$\QQ$}
\end{picture}

\noindent A very important feature here is that $R_*$ is
very small, and a powerful lever arises.
The measure $\QQ$ admits a simple theoretical representation
and an efficient empirical estimation procedure.
With this measure, we thus have a powerful tool
for the analysis of the relationship between $\PP$ and $\RR$.

Having estimated $\QQ$, we can move along the line
$$
\QQ+\RR\longrightarrow\PP.
$$
This enables us to estimate the expectations of various
variables with respect to~$\PP$, i.e. rewards.
The importance of this descends from the fact that the
direct empirical estimation of rewards is known to be
virtually impossible (see the discussion in~\cite{B95}
and the 20s example in~\cite{JPR05}).
The reason is that the expected returns are very small,
so that a slight misspecification of the data leads to
a significant \textit{relative} change in the estimated
expected returns.

According to~\eqref{i2}, the expected discounted profit
earned by the $i$-th asset is
$$
\EE_\PP\De S_1^i=-R_*\EE_\QQ\De S_1^i.
$$
Here $\De S_1^i=S_1^i-S_0^i$ is the discounted P\&L
produced by the $i$-th asset over the unit time period.
This is, in fact, the coherent-based analog of the SML relation.
The value $R_*$ is very small, so that $\EE_\QQ\De S_1^i$
is a medium size number, and it admits effective empirical
estimation procedures (the above lever at work!).
It is hard to estimate $R_*$ in practice, but for
the decision making purposes one typically needs to know
the values $\EE_\PP\De S_1^i$ up to multiplication by a joint
constant. Thus, it is sufficient to know only the values
$\EE_\QQ\De S_1^i$. For these values we provide simple
empirical estimation procedures.

Another way that we can follow is to move along the line
$$
\PP+\QQ\longrightarrow\RR.
$$
This methodology leads to the following coherent variant
of the empirical asset pricing:
$$
V=\frac{1}{1+R_*}\,\EE_\PP F+\frac{R_*}{1+R_*}\,\EE_\QQ F.
$$
Here $F$ is the payoff of a contingent claim and $V$ is
the price of~$F$.
A similar representation is provided for sensitivity
coefficients.
The expressions for~$V$ and for the sensitivity
coefficients admit efficient empirical estimation
procedures.

\skb
\textbf{3. Structure of the paper.}
In Section~\ref{CRM}, we recall basic definitions and
facts related to coherent risk measures.

In Section~\ref{CAPM}, we first consider the Markowitz-type
optimization problem with risk measured in a coherent way
and discuss the form of the efficient frontier.
Then we present the coherent-based variant of the CAPM.

In Section~\ref{REAP}, we introduce the extreme measure
of the market, describe its relationship with the
No Better Choice pricing, and apply it to estimating
rewards as well as to the empirical asset pricing.

Section~\ref{SC} contains the conclusions.

\section{Coherent Risk Measures}
\label{CRM}

\textbf{1. Coherent risk measures.}
Let $(\Omega,\F,\PP)$ be a probability space.
According to the definition introduced in~\cite{ADEH97},
\cite{ADEH99}, a \textit{coherent risk measure} is a
map $\rho\colon L^\infty\to\R$ (recall that $L^\infty$
is the space of bounded random variables)
satisfying the properties:
\begin{mitemize}
\item[(a)] (Subadditivity) $\rho(X+Y)\le\rho(X)+\rho(Y)$;
\item[(b)] (Monotonicity) If $X\le Y$, then
$\rho(X)\ge\rho(Y)$;
\item[(c)] (Positive homogeneity) $\rho(\la X)=\la\rho(X)$ for
$\la\in\R_+$;
\item[(d)] (Translation invariance)
$\rho(X+m)=\rho(X)-m$ for $m\in\R$;
\item[(e)] (Fatou property) If $|X_n|\le1$,
$X_n\xra{\PP}X$, then $\rho(X)\le\liminf_n\rho(X_n)$.
\end{mitemize}

The representation theorem proved in~\cite{ADEH97},
\cite{ADEH99} for the case of a finite~$\Omega$ and
in~\cite{D02} for the general case states that $\rho$
is a coherent risk measure if and only if there exists
a non-empty set $\DDD$ of probability measures absolutely
continuous with respect to~$\PP$ such that
\begin{equation}
\label{crm1}
\rho(X)=-\inf_{\QQ\in\DDD}\EE_\QQ X,\quad X\in L^\infty.
\end{equation}

When dealing with theory, we need to define coherent
risks not only on bounded random variables, but on
unbounded ones as well (most distributions used in
theory, like the normal or the lognormal ones, are unbounded).
For this, we take representation~\eqref{crm1} as the
definition, i.e. we define a coherent risk measure~$\rho$
on the space $L^0$ of all random variables as a map
\begin{equation}
\label{crm2}
\rho(X):=-\inf_{\QQ\in\DDD}\EE_\QQ X,\quad X\in L^0,
\end{equation}
where $\DDD$ is a set of probability measures absolutely
continuous with respect to~$\PP$ (this way to extend
coherent risks to~$L^0$ was proposed in~\cite{C061}).\footnote{The
expectation $\EE_\QQ X$ is understood here in the
generalized sense: $\EE_\QQ X:=\EE_\QQ X^+-\EE_\QQ X^-$,
where $X^+=\max\{X,0\}$, $X^-=\max\{-X,0\}$.
Thus, $\EE_\QQ$ and $\rho$ take on values in
$[-\infty,\infty]$.}

A set $\DDD$, for which~\eqref{crm2} is true, is not
unique (for example, $\DDD$ and its convex hull define
the same risk measure).
However, there exists the largest such set.
It consists of the measures~$\QQ$ absolutely continuous
with respect to~$\PP$ such that $\EE_\QQ X\ge-\rho(X)$
for any~$X$. We call it the \textit{determining set}
of~$\rho$.

For more information on coherent risk measures, we refer
to~\cite{CM061}, \cite{D05}, and~\cite[Ch.~4]{FS04}.

\skb
\textbf{2. Examples.}
Let us give examples of four most natural classes
of coherent risk measures: Tail V@R, Weighted V@R,
Beta V@R, and Alpha V@R.

\textit{Tail V@R of order} $\la\in(0,1]$
(the terms \textit{Average V@R},
\textit{Conditional V@R}, \textit{Expected Shortfall},
and \textit{Expected Tail Loss} are also used)
is the coherent risk measure~$\rho_\la$ corresponding
to the determining set
$$
\DDD_\la=\Bigl\{\QQ:\frac{d\QQ}{d\PP}
\le\la^{-1}\Bigr\}.
$$
If $X$ has a continuous distribution, then
$$
\rho_\la(X)=-\EE(X\cond X\le q_\la(X)),
$$
where $q_\la(X)$ is the $\la$-quantile of~$X$.
This motivates the term Tail V@R.
For a detailed study of this risk measure, we refer
to~\cite{AT02}, \cite[Sect.~2]{CM061},
and \cite[Sect.~4.4]{FS04}.

\textit{Weighted V@R with the weighting measure~$\mu$}
(the term \textit{spectral risk measure} is also used),
where $\mu$ is a probability measure on $(0,1]$,
is the coherent risk measure~$\rho_\mu$ defined~as
$$
\rho_\mu(X)=\int_{(0,1]}\rho_\la(X)\mu(d\la).
$$
One can check that this is indeed a coherent risk measure.
Its determining set will be denoted by~$\DDD_\mu$.
Weighted V@R admits several equivalent representations.
One of the most convenient representations is:
$$
\rho_\mu(X)=-\int_0^1 q_x(X)\psi_\mu(x)dx,
$$
where
\begin{equation}
\label{crm3}
\psi_\mu(x)=\int_{[x,1]}\la^{-1}\mu(d\la),
\quad x\in[0,1].
\end{equation}
In particular, let $\Omega=\{1,\dots,T\}$ and
$X(t)=x_t$. Let $x_{(1)},\dots,x_{(T)}$ be the values
$x_1,\dots,x_T$ in the increasing order.
Define $n(t)$ through the equality $x_{(t)}=x_{n(t)}$.
Then
$$
\rho_\mu(X)=-\sum_{t=1}^T x_{n(t)}
\int_{z_{t-1}}^{z_t}\psi_\mu(x)dx,
$$
where $z_t=\sum_{i=1}^t\PP\{n(i)\}$.
This formula provides a simple empirical estimation
procedure for~$\rho_\mu$.
For a detailed study of Weighted V@R, we refer
to~\cite{A02}, \cite{A04}, \cite{C05e},
\cite[Sect.~2]{CM061}, and~\cite[Sect.~4.6, 4.7]{FS04}.

\textit{Beta V@R with parameters}
$\al\in(-1,\infty)$, $\be\in(-1,\al)$ is the
Weighted V@R with the weighting measure
$$
\mu_{\al,\be}(dx)={\rm B}(\be+1,\al-\be)^{-1}x^\be(1-x)^{\al-\be-1}dx,
\quad x\in[0,1].
$$
As shown in~\cite{CM061}, for $\al,\be\in\N$, Beta V@R
admits the following simple representation
$$
\rho_{\al,\be}(X)=
-\EE\Bigl[\frac{1}{\be}\sum_{i=1}^\be X_{(i)}\Bigr],
$$
where $X_{(1)},\dots,X_{(\al)}$ are the order statistics
obtained from independent copies $X_1,\dots,X_\al$ of~$X$.
This representation provides a very convenient way
for the empirical estimation of~$\rho_{\al,\be}$.
For a detailed study of this risk measure,
see~\cite[Sect.~2]{CM061}.

\textit{Alpha V@R} is obtained from Beta V@R by fixing
$\be=1$. Clearly, if $\al\in\N$, then
$$
\rho_\al(X)=-\EE\min_{i=1,\dots,\al}X_i,
$$
where $X_1,\dots,X_\al$ are independent copies of~$X$.

In our opinion, the most important classes of coherent
risk measures are: Alpha V@R, Beta V@R, and Weighted V@R.

\skb
\textbf{3. $L^1$-spaces.}
For technical purposes, we need to recall the
definition of the \textit{strong $L^1$-space} associated
with a coherent risk measure~$\rho$:
$$
L_s^1(\DDD)=\Bigl\{X\in L^0:\lim_{n\to\infty}
\sup_{\QQ\in\DDD}\EE_\QQ|X|I(|X|>n)=0\Bigr\},
$$
where $\DDD$ is the determining set of~$\rho$.

Let us provide two examples.
For Weighted V@R,
$$
L_s^1(\DDD_\mu)
=\{X\in L^0:\rho_\mu(X)<\infty,\,\rho_\mu(-X)<\infty\}
$$
(see~\cite[Subsect.~2.2]{C061}).
The right-hand side of this equality was called
in~\cite{C061} the \textit{weak $L^1$-space}.
It has a clear financial interpretation: this is the
set of random variables such that their risk is finite
and the risk of their negatives is finite.

For Beta V@R with $\be>0$
(in particular, for Alpha V@R), $L_s^1$ has a very simple
form: it coincides with the space $L^1$ of $\PP$-integrable
random variables (see~\cite[Sect.~3]{CM061}).

\skb
\textbf{4. Extreme measures and risk contributions.}
Let $\rho$ be a coherent risk measure with the determining
set~$\DDD$ and $W\in L_s^1(\DDD)$.
From the financial point of view, $W$ is the
discounted P\&L produced
by some portfolio (for example, the portfolio of some
firm) over the unit time period.
A measure~$\QQ\in\DDD$ is called an \textit{extreme measure}
for~$W$ if $\EE_\QQ W=-\rho(W)$.
The set of extreme measures will be denoted by $\EX_\DDD(W)$.

As shown in~\cite[Subsect.~2.3]{C061}, $\EX_\DDD(W)$ is
non-empty provided that $\DDD$ is $L^1$-closed and
uniformly integrable.
The latter condition is automatically satisfied by~$\DDD_\mu$
(this follows from the explicit representations of~$\DDD_\mu$
provided in~\cite{CD03}, \cite[Sect.~4]{C05e}).

Let us give two examples.
If $W\in L_s^1(\DDD_\mu)$ has a
continuous distribution, then $\EX_{\DDD_\mu}(W)$
consists of a unique measure
\begin{equation}
\label{crm4}
\QQ_\mu(W)=\psi_\mu(F(W))\PP,
\end{equation}
where $\psi_\mu$ is given by~\eqref{crm3} and $F$ is the
distribution function of~$W$ (for the proof,
see~\cite[Sect.~6]{C05e}).

Let $\Omega=\{1,\dots,T\}$ and $W(t)=w_t$.
Assume that all the values $w_t$ are different.
Let $w_{(1)},\dots,w_{(T)}$ be these values
in the increasing order. Define $n(t)$ through the equality
$w_{(i)}=w_{n(i)}$. Then $\EX_{\DDD_\mu}(W)$ consists
of a unique measure $\QQ_\mu(W)$ given by
\begin{equation}
\label{crm5}
\QQ_\mu(W)\{n(t)\}=\int_{z_{t-1}}^{z_t}\psi_\mu(x)dx,
\end{equation}
where $\psi_\mu$ is defined by~\eqref{crm3} and
$z_t=\sum_{i=1}^t\PP\{n(i)\}$ (for the proof,
see~\cite[Sect.~5]{CM061}).

\skm
The notion of an extreme measure is closely connected
with the notion of \textit{risk contribution}.
Let $X\in L^0$ be the discounted P\&L produced by some
trade over the unit time period.
According to~\cite{C061}, the risk contribution of~$X$
to~$W$ is defined as
$$
\rho^c(X;W)=-\inf_{\QQ\in\EX_\DDD(W)}\EE_\QQ X.
$$
As $\EX_\DDD(W)$ is typically a singleton (see the examples
above), $\rho^c(X;W)$ is typically linear in~$X$.
The relevance of the above definition is seen from the
following result: if $X\in L_s^1(\DDD)$, then
\begin{equation}
\label{crm6}
\rho^c(X;W)=\lim_{\eps\da0}\eps^{-1}(\rho(W+\eps X)-\rho(W))
\end{equation}
(for the proof, see~\cite[Subsect.~2.5]{C061}).

Let us provide four examples.
If $W\in L_s^1(\DDD_\mu)$ has a
continuous distribution, then, according to~\eqref{crm4},
$$
\rho_\mu^c(X;W)
=-\EE_{\QQ_\mu(W)}X
=-\EE_\PP\psi_\mu(F(W))X.
$$

Let $\Omega=\{1,\dots,T\}$, $X(t)=x_t$, $W(t)=w_t$.
Assume that all the values $w_t$ are different.
Let $w_{(1)},\dots,w_{(T)}$ be these values
in the increasing order.
Define $n(t)$ through the equality $w_{(t)}=w_{n(t)}$.
Then, according to~\eqref{crm5},
\begin{equation}
\label{crm7}
\rho_\mu^c(X;W)
=-\EE_{\QQ_\mu(W)}X
=-\sum_{t=1}^T x_{n(t)}\int_{z_{t-1}}^{z_t}\psi_\mu(x)dx,
\end{equation}
where $\psi_\mu$ is given by~\eqref{crm3} and
$z_t=\sum_{i=1}^t\PP\{n(i)\}$.

Let $\al,\be\in\N$, $X,W\in L^1$, and
suppose that $W$ has a continuous distribution.
Let $(X_1,W_1),\dots,(X_\al,W_\al)$ be independent copies
of $(X,W)$ and let $W_{(1)},\dots,W_{(T)}$ be the corresponding
order statistics (i.e. the values $W_1,\dots,W_T$ in the
increasing order).
Define random variables $n(i)$ through the equality
$W_{(i)}=W_{n(i)}$ (as $W$ has a continuous distribution,
all the values $W_1,\dots,W_\al$ are a.s. different, so
that $n(i)$ is a.s. determined uniquely).
Then, as shown in~\cite[Sect.~5]{CM061},
\begin{equation}
\label{crm8}
\rho_{\al,\be}^c(X;W)=-\EE\Bigl[\frac{1}{\be}\sum_{i=1}^\be
X_{n(i)}\Bigr].
\end{equation}

Let $\al\in\N$,  $X,W\in L^1$, and
suppose that $W$ has a continuous distribution.
It follows from the above example that
\begin{equation}
\label{crm9}
\rho_\al^c(X;W)=-\EE X_{\argmin\limits_{i=1,\dots,\al}W_i},
\end{equation}
where $(X_1,W_1),\dots,(X_\al,W_\al)$ are independent
copies of $(X,W)$.

For more information on extreme measures and risk
contributions, we refer to~\cite{C061}, \cite{C062},
and~\cite{CM061}.

\section{CAPM}
\label{CAPM}

\textbf{1. Efficient frontier.}
Let $(\Omega,\F,\PP)$ be a probability space.
Let $\rho$ be a coherent risk measure with the determining
set~$\DDD$.
Let $S_0=(S_0^1,\dots,S_0^d)$ be the vector of prices of
assets $1,\dots,d$ at time~0 and $S_1=(S_1^1,\dots,S_1^d)$
be the vector of their discounted prices at time~1.
We will assume that $S_1^i\in L_s^1(\DDD)\cap L^1$.
An investor's strategy is described by a vector
$h=(h^1,\dots,h^d)$, whose $i$-th component means the
amount of assets of the type~$i$ bought by the investor
($h^i$ might be negative, which corresponds to the short
selling).
The discounted P\&L produced by a strategy $h$ is
$\sum_i h^i(S_1^i-S_0^i)=\lb h,\De S_1\rb$.

Let us consider the following Markowitz-type optimization
problem:
\begin{equation}
\label{capm1}
\begin{cases}
\EE_\PP\lb h,\De S_1\rb\longrightarrow\max,\\
h\in\R^d,\\
\rho(\lb h,\De S_1\rb)\le c,
\end{cases}
\end{equation}
where $c\in(0,\infty)$.
We are dealing with the discounted P\&Ls rather than with
returns because an agent need not really invest the money in order
to take a long/short position in an asset; he/she can
either borrow the money by posing a collateral or take
a long/short position in the futures again by posing a
collateral. We will assume that the collateral required
to support a strategy~$h$ is dominated by
$\rho(\lb h,\De S_1\rb)$.
Thus, if an agent possesses a capital $w$, he/she
divides it into two parts: a sum $c$ for the risky
investments and a sum $w-c$ for the risk-free investment.
Then he/she chooses a strategy~$h$ such that
$\rho(\lb h,\De S_1\rb)\le c$ and, using this sum as
a collateral, borrows the money needed to buy/sell short
the corresponding assets. Thus, the sum invested
by him/her into risky assets is several times larger
than~$c$.

Problems of type~\eqref{capm1} were considered
in~\cite{A04}, \cite{C062}, \cite{RU00}, \cite{RUZ05}.
Let us describe a geometric solution proposed
in~\cite[Subsect.~2.2]{C062}.
We will assume that $\rho(\lb h,\De S_1\rb)>0$ for any
$h\in\R^d\setminus\{0\}$, which means that any trade has
a strictly positive risk. Consider the set
$$
G=\cl\{\EE_\QQ\De S_1:\QQ\in\DDD\},
$$
where ``$\cl$'' denotes the closure.
This set is a convex compact in~$\R^d$ containing 0 as
an inner point.
In~\cite{C061}, it was termed the \textit{generator}
of~$\De S_1$ and~$\rho$. Its role is seen from the line
$$
\rho(\lb h,\De S_1\rb)
=-\inf_{\QQ\in\DDD}\EE_\QQ\lb h,\De S_1\rb
=-\inf_{\QQ\in\DDD}\lb h,\EE_\QQ\De S_1\rb
=-\min_{x\in G}\lb h,x\rb,\quad h\in\R^d.
$$
Let $T$ be the intersection of the ray
$(\EE_\PP\De S_1,0)$ with the border of~$G$.
Let $N$ be the set of the inner normals to~$G$ at the
point~$T$. As shown in~\cite{C062}, the set of solutions
of~\eqref{capm1}~is
$$
\bigl\{h\in N:\min_{x\in G}\lb h,x\rb=-c\bigr\}.
$$
This set is non-empty.
In general, the border of $G$ might have a break at the
point~$T$, so that this set may contain more than one
strategy (for instance, this might happen if one of the
assets $1,\dots,d$ is an option written on another;
see~\cite[Subsect.~2.2]{C062} for details).
However, if $S_1^1,\dots,S_1^d$ have a joint density,
$\rho=\rho_\mu$, and the support of~$\mu$ is the whole
interval $[0,1]$, then the optimal strategy~$h_*$
is unique (see~\cite[Sect.~5]{C05e}).
In this case the efficient frontier on the $(E,\rho)$-plane
is $(R_*c,c:c\in\R_+)$, where
$$
R_*=\frac{\EE_\PP\lb h_*,\De S_1\rb}{\rho(\lb h_*,\De S_1\rb)}.
$$

\begin{figure}[!h]
\begin{picture}(150,52.5)(-80,-35)
\put(-28.5,-16){\includegraphics{capm.1}}
\put(-22,-14){\small $T$}
\put(0.2,-3.2){\small $0$}
\put(4,0.5){\small $\EE_\PP\De S_1$}
\put(-14,0){\small $h_*$}
\put(-1.5,7.7){\small $G$}
\put(-50,-30){\small\textbf{Figure~1.}
Geometric solution of the optimization problem}
\end{picture}
\end{figure}

\begin{figure}[!h]
\begin{picture}(150,42.5)(-45,-15)
\put(0,-0.3){\includegraphics{capm.2}}
\put(0,0){\vector(1,0){70}}
\put(0,0){\vector(0,1){25}}
\put(69,-3){$\rho$}
\put(-4.5,23){\small $E$}
\put(35.5,2){\scalebox{0.75}{Attainable portfolios}}
\put(25,5){\rotatebox{10}{\small$E=R_*\rho$}}
\put(8,-12){\small\textbf{Figure~2.}
The efficient frontier}
\end{picture}
\end{figure}

The above geometric solution provides an insight into
the form of the optimal portfolio and is convenient if we
have a model for the joint distribution of
$S_1^1,\dots,S_1^d$. For example, if this vector is
Gaussian with a non-degenerate covariance
matrix~$C$, then $h_*$ coincides up to a constant with
$C^{-1}\EE_\PP\De S_1$ (see~\cite[Subsect.~2.2]{C062}).
However, if we have no model for $S_1^1,\dots,S_1^d$, but
rather want to solve~\eqref{capm1} empirically, then we
can rewrite it in the equivalent form
$$
\begin{cases}
\rho(\lb h,\De S_1\rb)\longrightarrow\min,\\
h\in\R^d,\\
\EE_\PP\lb h,\De S_1\rb=1.
\end{cases}
$$
This is the problem of minimizing a convex functional
over an affine space. For example, if $\rho$ is Alpha V@R,
Beta V@R, or Weighted V@R, we can use the empirical
estimation procedures for $\rho(\lb h,\De S_1\rb)$
described in Section~\ref{CRM} and approach the above
problem numerically.

\skb
\textbf{2. Security market line.}
Suppose that the optimal strategy~$h_*$ for~\eqref{capm1}
is unique. We will assume that the set of extreme
measures $\EX_\DDD(\lb h_*,\De S_1\rb)$ consists of a
unique measure~$\QQ$.
As shown by~\eqref{crm4},
this assumption is satisfied if $\rho=\rho_\mu$ and
$S_1^1,\dots,S_1^d$ have a joint density.
Fix $i\in\{1,\dots,d\}$ and set
$h(\eps)=(h_*^1,\dots,h_*^{i-1},h_*^i+\eps,
h_*^{i+1},\dots,h_*^d)$. According to~\eqref{crm6},
$$
\frac{d}{d\eps}\Bigl|_{\eps=0}\rho(\lb h(\eps),\De S_1\rb)
=-\EE_\QQ\De S_1^i.
$$
Obviously,
$$
\frac{d}{d\eps}\Bigl|_{\eps=0}\EE_\PP\lb h(\eps),\De S_1\rb
=\EE_\PP\De S_1^i.
$$
Since $h_*$ delivers the maximum of
$\EE_\PP\lb h,\De S_1\rb/\rho(\lb h,\De S_1\rb)$ and
this maximum equals~$R_*$, we get
\begin{equation}
\label{capm3}
\EE_\PP\De S_1^i=-R_*\EE_\QQ\De S_1^i,\quad i=1,\dots,d.
\end{equation}
Consider the returns
\begin{align*}
r^i&=\frac{(1+r_f)S_1^i-S_0^i}{S_0^i},\quad i=1,\dots,d,\\
r_*&=\frac{(1+r_f)\lb h_*,S_1\rb-\lb h_*,S_0\rb}{\lb h_*,S_0\rb},
\end{align*}
where $r_f$ is the risk-free interest rate (recall that
$S_1^1,\dots,S_1^d$ are the discounted prices).
In view of the equalities
\begin{align*}
r^i-r_f&=\frac{1+r_f}{S_0^i}\,\De S_1^i,\quad i=1,\dots,d,\\
r_*-r_f&=\frac{1+r_f}{\lb h_*,S_0\rb}\,\lb h_*,\De S_1\rb,
\end{align*}
we see that~\eqref{capm3} transforms into
$$
\EE_\PP(r^i-r_f)
=\be^i\EE_\PP(r_*-r_f),\quad i=1,\dots,d,
$$
where
$$
\be^i=\frac{\EE_\QQ(r^i-r_f)}{\EE_\QQ(r_*-r_f)},
\quad i=1,\dots,d.
$$
This is the analog of the SML relation.
The coefficients $\be^i$ admit simple empirical estimation
procedures described in Section~\ref{REAP}.

\skb
\textbf{3. Equilibrium.}
Suppose now that there are $N$ agents in the economy.
The $n$-th agent has an initial endowment $W_n\in(0,\infty)$
and a risk aversion coefficient $a_n\in(0,\infty)$.
Suppose that the $n$-th agent has preferences that are
linear in reward and quadratic in risk,\footnote{This kind of
assumption is typical for the economic theory.
For example, if the agent is using expected utility
with the utility function $U(x)=A-e^{-\la x}$,
then the certainty equivalent of a Gaussian random
variable~$X$ with mean~$a$ and variance~$\si^2$ is
$a-\frac{\la}{2}\si^2$. Furthermore, if $\rho(X)$
depends only on the distribution of~$X$
(for example, this assumption is satisfied by Weighted V@R),
then there exists a constant $\ga>0$ such
that, for a Gaussian random variable~$X$ with
variance~$\si^2$, $\rho(\wl X)=\ga\si$,
where $\wl X=X-\EE X$ is the centered version of~$X$.
Let us finally remark that in the third line of~\eqref{capm1}
we can use the centered versions of $\lb h,\De S_1\rb$
without essentially changing the problem because
$\EE_\PP\De S_1$ is close to~0.} i.e. he/she is solving the problem
\begin{equation}
\label{capm4}
\begin{cases}
\ds\EE_\PP\lb h,\De S_1\rb-\frac{a_n}{2W_n}\,
\rho^2(\lb h,\De S_1\rb)\longrightarrow\max,\\
h\in\R^d.
\end{cases}
\end{equation}
Let $h_*$ be the solution of~\eqref{capm1} with $c=1$
(we assume that it is unique).
Clearly, the solution of~\eqref{capm4} is given by
$h_n=c_n h_*$ with some positive constant $c_n$.
This constant is found by solving the problem
$$
c R_*-c^2\frac{a_n}{2W_n}\longrightarrow\max,
$$
which yields $c_n=R_*W_n/a_n$.

In order to find $R_*$ using the equilibrium considerations,
we will consider $S_1^1,\dots,S_1^d$ and the market
portfolio as known variables and $S_0^1,\dots,S_0^d$ as
unknown variables (derived from the equilibrium).
The market portfolio has the form
$$
H_*=\sum_{n=1}^N c_n h_*=Ch_*.
$$
In particular, the portfolio of the $n$-th agent is
$$
h_n=H_*\frac{W_n}{a_n}\biggl(\sum_{n=1}^N\frac{W_n}{a_n}\biggr)^{-1}.
$$
Using the equality $\rho(\lb h_*,\De S_1\rb)=1$, we can
write
$$
R_*
=\frac{\EE_\PP\lb H_*,\De S_1\rb}{\rho(\lb H_*,\De S_1\rb)}
=\frac{\EE_\PP\lb H_*,\De S_1\rb+\rho(\lb H_*,\De S_1\rb)}
{\rho(\lb H_*,\De S_1\rb)}-1
=\frac{\EE_\PP\lb H_*,S_1\rb+\rho(\lb H_*,S_1\rb)}{C}-1.
$$
On the other hand,
$$
C=\sum_{n=1}^N c_n=R_*\sum_{n=1}^N\frac{W_n}{a_n}.
$$
From the last two equalities we get the equation for~$R_*$:
$$
R_*^2+R_*-\biggl(\sum_{n=1}^N\frac{W_n}{a_n}\biggr)^{-1}
[\EE_\PP\lb H_*,S_1\rb+\rho(\lb H_*,S_1\rb)]=0,
$$
which yields
$$
R_*=-\frac12+\frac12\biggl(1+4\biggl(\sum_{n=1}^N\frac{W_n}{a_n}\biggr)^{-1}
[\EE_\PP\lb H_*,S_1\rb+\rho(\lb H_*,S_1\rb)]\biggr)^{\frac12}.
$$

\section{Rewards and Empirical Asset Pricing}
\label{REAP}

\textbf{1. Extreme measure and contact measure.}
Consider the framework of the previous section and assume
that $\DDD$ is $L^1$-closed and uniformly integrable
(as mentioned above, this assumption is satisfied, for
example, by~$\DDD_\mu$).
Define by~$\RRR$ the set of the \textit{risk-neutral
measures}, i.e. the measures~$\RR$ absolutely continuous
with respect to~$\PP$ such that $\EE_\RR\De S_1=0$.
It follows from the results of~\cite[Subsect.~3.2]{C061}
that
$$
\biggl(\frac{1}{1+R}\,\{\PP\}
+\frac{R}{1+R}\,\DDD\biggr)\cap\RRR=\emp
$$
for $R<R_*$ and
$$
\biggl(\frac{1}{1+R_*}\,\{\PP\}
+\frac{R_*}{1+R_*}\,\DDD\biggr)\cap\RRR\ne\emp.
$$
We will call the latter set the set of \textit{contact measures}.

Let $\QQ\in\EX_\DDD(\lb h_*,\De S_1\rb)$.
Consider the measure
$$
\RR=\frac{1}{1+R_*}\,\PP+\frac{R_*}{1+R_*}\,\QQ.
$$
It follows from~\eqref{capm3} that $\RR\in\RRR$, so that
$\RR$ is a contact measure.
Conversely, if $\QQ$ is arbitrary and the measure~$\RR$
is a contact measure, then $\QQ\in\DDD$ and
$$
\EE_\PP\lb h_*,\De S_1\rb=-R_*\EE_\QQ\lb h_*,\De S_1\rb.
$$
On the other hand,
$$
\EE_\PP\lb h_*,\De S_1\rb=R_*\rho(\lb h_*,\De S_1\rb).
$$
These equalities combined together show that
$\QQ\in\EX_\DDD(\lb h_*,\De S_1\rb)$.
Thus, we have proved the equality
$$
\frac{1}{1+R_*}\,\{\PP\}+\frac{R_*}{1+R_*}\,\CCC
=\EX_\DDD(\lb h_*,\De S_1\rb),
$$
where $\CCC$ is the set of contact measures.
If we assume as in the previous section that
$\EX_\DDD(\lb h_*,\De S_1\rb)$ consists of a unique
measure~$\QQ$, then we obtain that the set of contact
measures consists of a unique measure~$\RR$, and the
following relation holds:
$$
\RR=\frac{1}{1+R_*}\,\PP+\frac{R_*}{1+R_*}\,\QQ.
$$

\begin{figure}[!h]
\begin{picture}(150,65)(-77.5,-37.5)
\put(-30.2,-25.3){\includegraphics{capm.3}}
\multiput(0,0)(1,0){30}{\line(1,0){0.5}}
\put(-1,-4){\small $\PP$}
\put(2,-4){\small $\RR$}
\put(28,-4){\small $\QQ$}
\put(-1.8,9){\small $\DDD$}
\put(28.2,17){\small $\RRR$}
\put(-19,-20){\scalebox{0.8}{$\frac{1}{1+R_*}\,\{\PP\}
+\frac{R_*}{1+R_*}\,\DDD$}}
\put(-36,-35){\small\textbf{Figure 3.}
The joint arrangement of $\PP$, $\QQ$, $\RR$}
\end{picture}
\end{figure}

The notion of a contact measure is closely connected
with the \textit{No Better Choice} (\textit{NBC})
pricing technique introduced in~\cite[Subsect.~3.1]{C062}.
Let $F$ be a random variable meaning the discounted
payoff of some contingent claim.
A real number~$x$ is an \textit{NBC price} of~$F$~if
$$
\sup_{h\in\R^d\!,\,g\in\R}\frac{\EE_\PP(\lb h,\De S_1\rb+g(F-x))}%
{\rho(\lb h,\De S_1\rb+g(F-x))}
=\sup_{h\in\R^d}\frac{\EE_\PP\lb h,\De S_1\rb}{\rho(\lb h,\De S_1\rb)}.
$$
In other words, $x$ is an NBC price if the
incorporation of the $d+1$-st asset with the initial
price~$x$ and the terminal
price~$F$ does not increase the optimal value in
problem~\eqref{capm1}.
As shown in~\cite{C062}, the set of NBC prices coincides
with the interval $\{\EE_\QQ F:\QQ\in\CCC\}$.
In our situation $\CCC=\{\RR\}$, so that the NBC price
is $\EE_\RR F$.

\skb
\textbf{2. Rewards.}
Formula~\eqref{capm3} provides a convenient tool for
estimating rewards.
The measure~$\QQ$ is the extreme measure for
$\lb H_*,\De S_1\rb$, where $H_*$ is the market portfolio
(note that if two random variables coincide up to
multiplication by a positive
constant, their extreme measures coincide).
The random variable $\lb H_*,\De S_1\rb$ is the increment
of the overall market wealth over the unit time period.
Thus, $\EE_\QQ\De S_1^i$ is nothing but the risk
contribution of the $i$-th asset to the whole market.
So, \eqref{capm3} has the following meaning:
$$
\text{Reward}=\text{Price of risk}\times\text{Risk
contribution to the market}.
$$
It is typically hard to estimate~$R_*$, but for most
decision making purposes it is sufficient to know the
values $\EE_\PP\De S_1^i$ up to multiplication by a
positive constant. Thus, it is sufficient to estimate
$\EE_\QQ\De S_1^i$.

A good approximation to $\lb H_*,\De S_1\rb$
is the increment~$\De I_1=I_1-I_0$
of the market index like S\&P~500 over the unit time period.
If $\QQ$ is the extreme measure of $\De I_1$, then
$\EE_\QQ\De S_1^i$ is the risk contribution
of~$\De S_1^i$ to~$\De I_1$.
We have
\begin{align*}
\De S_1^i&=\frac{S_0^i}{1+r_f}\,(r^i-r_f),\\
\De I_1&=\frac{I_0}{1+r_f}\,(r_I-r_f),
\end{align*}
where $r^i$ (resp., $r_I$) is the return on the $i$-th asset
(resp., the index) and $r_f$ is the risk-free rate.
Hence,
\begin{align*}
\EE_\QQ\De S_1^i
&=\frac{S_0^i}{1+r_f}\,\EE_\QQ(r^i-r_f)\\[1mm]
&=-\frac{S_0^i}{1+r_f}\,\rho^c(r^i-r_f;I_0(r_I-r_f))\\[1mm]
&=-\frac{S_0^i}{1+r_f}\,\rho^c(r^i;r_I)
-\frac{S_0^i r_f}{1+r_f}.
\end{align*}
Theoretical methods of calculating risk contributions
were described in Section~\ref{CRM}.
Here we will discuss its practical estimation.
Below $\delta$ denotes the length of the unit time interval.

Suppose first that $\rho$ is Weighted V@R.
In order to estimate $\rho_\mu^c(r^i;r_I)$,
one should first generate the data set
$(x_1,y_1),\dots,(x_T,y_T)$ for $(r^i,r_I)$
and a probability measure~$\nu$ on this set.
This can be done (in particular) by one of the following
techniques:
\begin{mitemize}
\item[\bf 1.] \textit{Historical simulation.}
The values $(x_1,y_1),\dots,(x_T,y_T)$ are $T$ recent
realizations of $(r^i,r_I)$; $\nu$ is uniform.
For example, if $\delta$ is one day, these are $T$
recent daily returns of the $i$-th asset and of the index.
\item[\bf 2.] \textit{Weighted historical simulation.}
The values $(x_1,y_1),\dots,(x_T,y_T)$ are $T$ recent
realizations of $(r^i,r_I)$, while $\nu$ is a
measure giving more mass to recent realizations.
For example, a natural choice for $\nu$ is the geometric
distribution with a parameter~$\la\in[0.95,0.99]$.
\item[\bf 3.] \textit{Bootstrapped historical simulation.}
We split the time axis into small intervals of length
$n^{-1}\delta$ and create each $(x_t,y_t)$ as
$$
x_t=\prod_{k=1}^n(1+r_k^i)-1,\qquad
y_t=\prod_{k=1}^n(1+r_{Ik})-1,
$$
where $r_k^i$ and $r_{Ik}$ are the returns of the $i$-th
asset and of the index over $n$ randomly
chosen small intervals; $\nu$ is uniform.
This method can be combined with the weighting method:
recent small intervals can be drawn with a higher
probability than older ones (in this case $\nu$ is
still uniform).
\item[\bf 4.] \textit{Monte Carlo simulation.}
The values $(x_t,y_t)$ are drawn at random from a
distribution of $(r^i,r_I)$ estimated from the historic
data; $\nu$ is uniform. For example, they might be drawn
from an ARCH or GARCH model.
\end{mitemize}
Let $y_{(1)},\dots,y_{(T)}$ be the values $y_1,\dots,y_T$
in the increasing order. Define $n(t)$ through the
equality $y_{(t)}=y_{n(t)}$.
According to~\eqref{crm7}, an
estimate of $\rho_\mu^c(r^i;r_I)$ is provided~by
\begin{equation}
\label{reap1}
\rho_{\text{e}}^c(r^i;r_I)
=-\sum_{t=1}^T x_{n(t)}
\int_{z_{t-1}}^{z_t}\psi_\mu(x)dx,
\end{equation}
where $z_t=\sum_{i=1}^t\nu\{n(i)\}$ and $\psi_\mu$ is
given by~\eqref{crm3}.

If $\rho$ is Beta V@R with $\al,\be\in\N$,
then instead of the procedure described above one can
use a faster Monte Carlo procedure.
One should fix a number of trials $K\in\N$ and
generate independent draws
$((x_{kl},y_{kl});k=1,\dots,K,\,l=1,\dots,\al)$ of
$(r^i,r_I)$. This can be done using the data
selection methods 1--4 described above.
For example, if we use the weighted historical simulation,
we are drawing $(x_{kl},y_{kl})$ from the recent
$T$ realizations of $(r^i,r_I)$
($T$ might be equal to~$\infty$)
in accordance with the measure~$\nu$.
Let $l_{k1},\dots,l_{k\be}$ be the numbers
$l\in\{1,\dots,\al\}$ such that the corresponding $y_{kl}$
stand at the first $\be$ places (in the increasing order)
among $y_{k1},\dots,y_{k\al}$.
According to~\eqref{crm8}, an
estimate of $\rho_{\al,\be}^c(r^i;r_I)$ is provided by
\begin{equation}
\label{reap2}
\rho_{\text{e}}^c(r^i;r_I)
=-\frac{1}{K\be}\sum_{k=1}^K
\sum_{i=1}^\be x_{kl_{ki}}.
\end{equation}

If $\rho$ is Alpha V@R with $\al\in\N$, one should
generate $x_{kl},y_{kl}$ similarly and calculate
the array
$$
l_k=\argmin_{l=1,\dots,\al}y_{kl},\quad k=1,\dots,K.
$$
According to~\eqref{crm9}, an
estimate of $\rho_\al^c(r^i;r_I)$ is provided by
\begin{equation}
\label{reap3}
\rho_{\text{e}}^c(r^i;r_I)
=-\frac{1}{K}\sum_{k=1}^K x_{kl_k}.
\end{equation}

An advantage of Alpha V@R and Beta V@R over Weighted V@R
is that for these risk measures the above described
empirical estimation procedure does not require the
ordering of $y_1,\dots,y_T$ (the ordering of this set
requires $T\log_2 T$ operations; this is a particularly
unpleasant number for $T=\infty$, which is a typical
value for~$T$ in the weighted historical simulation).
Of course, Weighted V@R is a wider class, but Beta V@R
is already rather a flexible family, and we believe that
one can confine himself/herself to this class (or just
to the class Alpha V@R).

Let us finally remark that reward (unlike risk) is
linear, so for the estimation of the reward of a large
portfolio it is sufficient to estimate the reward of each
of its components.

\skb
\textbf{3. Empirical asset pricing.}
Let $F=f(S_1)$ be the discounted payoff of some contingent
claim, where $S_1$ is the terminal value of the underlying
asset (this is not a discounted value).
The unit time interval here is of order of several months.
Let $r=S_0^{-1}(S_1-S_0)$ denote the return on the asset.
As mentioned above, the NBC price of $F$ is given by
$$
V
=\EE_\RR f(S_1)
=\EE_\RR f(S_0(1+r))
=\frac{1}{1+R_*}\,\EE_\PP f(S_0(1+r))
+\frac{R_*}{1+R_*}\,\EE_\QQ f(S_0(1+r)).
$$
The sensitivity of $V$ with respect to~$S_0$ is
$$
\frac{\partial V}{\partial S_0}
=\frac{1}{1+R_*}\,\EE_\PP(1+r)f'(S_0(1+r))
+\frac{R_*}{1+R_*}\,\EE_\QQ(1+r)f'(S_0(1+r)).
$$

The empirical estimation procedures for these quantities
are similar to the procedures described above.
If $\rho$ is Weighted V@R, one should generate a
data set $(x_1,y_1),\dots,(x_T,y_T)$ for $(r,r_I)$
and a measure~$\nu$ in the same way as above.
An estimate of~$V$ (the sensitivities are estimated
in a similar way) is provided by
$$
V_{\text{e}}
=\frac{1}{1+R_*}\sum_{t=1}^T f(S_0(1+x_t))\nu(\{t\})
+\frac{R_*}{1+R_*}\sum_{t=1}^T f(S_0(1+x_{n(t)}))
\int_{z_{t-1}}^{z_t}\psi_\mu(x)dx,
$$
where $n(t)$ is the same as in~\eqref{reap1}.

If $\rho$ is Beta V@R, one should generate
$((x_{kl},y_{kl});k=1,\dots,K,\,l=1,\dots,\al)$ in the
same way as above. An estimate of~$V$ is provided by
$$
V_{\text{e}}
=\frac{1}{K\al(1+R_*)}\sum_{k=1}^K\sum_{l=1}^\al
f(S_0(1+x_{kl}))
+\frac{R_*}{K\be(1+R_*)}\sum_{k=1}^K\sum_{i=1}^\be
f(S_0(1+x_{kl_{ki}})),
$$
where $l_{ki}$ is the same as in~\eqref{reap2}.

If $\rho$ is Alpha V@R, then an empirical estimate
of~$V$ is provided by
$$
V_{\text{e}}
=\frac{1}{K\al(1+R_*)}\sum_{k=1}^K\sum_{l=1}^\al
f(S_0(1+x_{kl}))
+\frac{R_*}{K(1+R_*)}\sum_{k=1}^K f(S_0(1+x_{kl_k})),
$$
where $l_k$ is the same as in~\eqref{reap3}.

Let us now compare the proposed variant of the empirical
asset pricing with the classical one.
In the classical approach,
$$
V=\EE_\PP\psi f(S_0(1+r)),
$$
where $\psi$ is the risk-aversion adjustment.
A typical choice is: $\psi=cU'(W_1)$, where
$W_1$ is the wealth of the economy at time~$1$,
$U$ is the utility function of a representative
investor, and $c$ is the normalizing constant chosen in
such a way that $\EE_\PP\psi=1$.
An approximation to $W_1$ is $W_0I_1/I_0$,
where $W_0$ is the wealth of the economy at time~0 and
$I_n$ is the value of the market index like S\&P~500 at time~$n$.
Thus, the risk-aversion adjustment in the pricing formula is
\begin{equation}
\label{reap4}
V-\EE_\PP f(S_0(1+r))=\EE_\PP(\psi-1)f(S_0(1+r)).
\end{equation}
In contrast, in our approach
\begin{equation}
\begin{split}
\label{reap5}
V-\EE_\PP f(S_0(1+r))
&=-\frac{R_*}{1+R_*}\,\EE_\PP f(S_0(1+r))
+\frac{R_*}{1+R_*}\,\EE_\QQ f(S_0(1+r))\\[1mm]
&\approx-R_*\EE_\PP f(S_0(1+r))+R_*\EE_\QQ f(S_0(1+r))\\[1mm]
&=R_*\EE_\PP(1-\phi)f(S_0(1+r)),
\end{split}
\end{equation}
where $\phi=\frac{d\QQ}{d\PP}$.
Note the difference between~\eqref{reap4} and~\eqref{reap5}:
$\psi$ is close to~$1$, so that $\psi-1$ is small;
on the other hand, $\phi$ is far from~1, but $R_*$ is
small, so that $R_*(1-\phi)$ is also small.

\section{Summary and Conclusion}
\label{SC}

\textbf{1. Security market line.} We prove that the
expectations of the discounted P\&Ls provided by
different assets have the form
\begin{equation}
\label{sc1}
\EE_\PP\De S_1^i=-R_*\EE_\QQ\De S_1^i,\quad i=1,\dots,d,
\end{equation}
where $R_*$ is the reward/risk ratio for the market
portfolio and $\QQ$ is the extreme measure of the market.
In terms of returns,
$$
\EE_\PP(r^i-r_f)=\be^i\EE_\PP(r_*-r_f),\quad i=1,\dots,d,
$$
where $r^i$ (resp., $r_*$) is the return on the $i$-th
asset (resp., the market portfolio), $r_f$ is the
risk-free rate, and
$$
\be^i=\frac{\EE_\QQ(r^i-r_f)}{\EE_\QQ(r_*-r_f)},
\quad i=1,\dots,d.
$$

\skb
\textbf{2. Equilibrium.} In the equilibrium,
$$
R_*=-\frac12+\frac12\biggl(1+4\biggl(\sum_{n=1}^N\frac{W_n}{a_n}\biggr)^{-1}
[\EE_\PP\lb H_*,S_1\rb+\rho(\lb H_*,S_1\rb)]\biggr)^{\frac12},
$$
where $W_n$ is the initial endowment of the $n$-th investor,
$a_n$ is his/her risk aversion coefficient,
$H_*$ is the market portfolio, and $S_1^i$ is the
discounted value of the $i$-th asset at time~1.
The portfolio of the $n$-th agent is
$$
h_n=H_*\frac{W_n}{a_n}\biggl(\sum_{n=1}^N\frac{W_n}{a_n}\biggr)^{-1}.
$$

\skb
\textbf{3. Extreme measure and contact measure.}
In addition to the real-world measure and the risk-neutral
measure, we introduce the third one: the extreme measure
of the market. The three measures are related by the
equality
$$
\RR=\frac{1}{1+R_*}\,\PP+\frac{R_*}{1+R_*}\,\QQ.
$$
In the case of an incomplete market, the measure $\RR$
is a particular representative of the set
of risk-neutral measures: the contact measure.
Using it as the pricing kernel corresponds to the
No Better Choice pricing.

The reason why such an important measure as~$\QQ$ emerges
is the very nature of coherent risk measures: they are
based on probabilistic scenarios and thus give rise to
very important probability kernels like the extreme
measure and the contact one.

\skb
\textbf{4. Rewards.}
Equality~\eqref{sc1} provides a convenient tool for
estimating rewards.
Assuming that $\QQ$ is the extreme measure of the index,
we get
$$
\EE_\QQ\De S_1^i
=-\frac{S_0^i}{1+r_f}\rho^c(r^i;r_I)
-\frac{S_0^i r_f}{1+r_f},
$$
where $r_I$ is the return on the index.
In order to estimate $\rho^c(r^i;r_I)$ for the case,
where $\rho$ is Alpha V@R with $\al\in\N$
or Beta V@R with $\al,\be\in\N$, one should fix $K\in\N$
and generate independent draws
$((x_{kl},y_{kl});k=1,\dots,K,\,l=1,\dots,\al)$ of
$(r^i,r_I)$.
This can be done by one of the following techniques:
\begin{mitemize}
\item historical simulation;
\item weighted historical simulation;
\item bootstrapped historical simulation;
\item Monte Carlo simulation.
\end{mitemize}

An estimate of $\rho_{\al,\be}^c(r^i;r_I)$ is provided by
$$
\rho_{\text{e}}^c(r^i;r_I)
=-\frac{1}{K\be}\sum_{k=1}^K
\sum_{i=1}^\be x_{kl_{ki}},
$$
where $l_{k1},\dots,l_{k\be}$ are the numbers
$l\in\{1,\dots,\al\}$ such that the corresponding $y_{kl}$
stand at the first $\be$ places (in the increasing order)
among $y_{k1},\dots,y_{k\al}$.

An estimate of $\rho_\al^c(r^i;r_I)$ is provided by
$$
\rho_{\text{e}}^c(r^i;r_I)
=-\frac{1}{K}\sum_{k=1}^K x_{kl_k},
$$
where
$$
l_k=\argmin_{l=1,\dots,\al}y_{kl},\quad k=1,\dots,K.
$$

Let us remark that
\begin{mitemize}
\item for the risk measurement purposes, it is
important to estimate the risk contributions to the
firm (as mentioned in the introduction);
\item for the reward measurement purposes, it is
important to estimate the risk contributions to the
whole market (as seen from the above considerations).
\end{mitemize}

\skb
\textbf{5. Empirical asset pricing.}
Let $F=f(S_1)$ be the discounted payoff of some contingent
claim. Its price and sensitivity are given by
\begin{align*}
V&=\frac{1}{1+R_*}\,\EE_\PP f(S_0(1+r))
+\frac{R_*}{1+R_*}\,\EE_\QQ f(S_0(1+r)),\\[2mm]
\frac{\partial V}{\partial S_0}&=
\frac{1}{1+R_*}\,\EE_\PP(1+r)f(S_0(1+r))
+\frac{R_*}{1+R_*}\,\EE_\QQ(1+r)f(S_0(1+r)),
\end{align*}
where $r=S_0^{-1}(S_1-S_0)$.
These values admit simple empirical estimation procedures
similar to those described above.

\skb
To conclude, we would like to draw the reader's attention
to the following fact. The classical CAPM is based on
the mean--variance analysis.
The classical empirical asset pricing is based on the
expected utility.
The classical risk measurement employs V@R.
In contrast, in our approach both the CAPM and the
empirical asset pricing are based on the coherent risk.
Of course, risk measurement can be based on the coherent risk.
Thus, a very big advantage of coherent risk measures is
their universality: they can be used
\begin{mitemize}
\item to measure risk;
\item for the pricing purposes (various theoretical
and empirical techniques are available);
\item for the decision making.
\end{mitemize}

\clearpage

\end{document}